\documentclass[12pt]{elsarticle}

\usepackage{color}

\begin{document}

\begin{frontmatter}

\title{On the Study of the General Group Classification of Systems of Linear Second-Order
Ordinary Differential Equations}

\author[SUN]{S.V. Meleshko}

\ead{sergey@math.sut.ac.th}

\author[DUT]{S. Moyo}

\ead{moyos@dut.ac.za}

\address[SUN]{Suranaree University of Technology, School of Mathematics, Nakhon
Ratchasima 30000, Thailand}

\address[DUT]{Durban University of Technology, Department of Mathematics \& Institute for Systems Science, P O Box 1334, Steve
Biko Campus, Durban 4000, South Africa}

\begin{abstract}
In this paper we study the general group classification of systems of linear second-order
ordinary differential equations inspired from earlier works and recent results on the group classification of such systems. Some interesting results and subsequent Theorem arising from this particular study are discussed here. This paper considers the study of irreducible systems of second-order ordinary differential equations.
\end{abstract}

\end{frontmatter}

\section{Introduction}

The appearance of systems of ordinary differential equations in the modeling of natural
phenomena has led to a vast interest in the study of their properties and theoretical aspects
that include their algebraic and symmetry properties. The existence of symmetries in a particular system leads to the possibility of reducing the order of the system or computing a general solution through quadratures.  This has been part of the reason for focusing on such systems in the current study. We study the general group classification of systems of linear second-order ordinary differential equations motivated by recent results obtained in \cite{bk:MoyoMeleshkoOguis[2013],bk:MoyoMeleshkoOguis[2014],%
bk:MkhizeMoyoMeleshko[2014],bk:SuksernMoyoMeleshko[2014]}.  Linear equations play an important role in many applications where they occur in a disguised form. When studying their symmetry properties it is always preferred to express them in their simplest equivalent form. It is important to note here that symmetry properties are invariant with respect to the change of the dependent and independent variables and hence are not affected by working with the equivalent form of a given system. The group classification problem involves classifying given differential equations with respect to arbitrary elements.  Here we use the algebraic algorithm in the group classification approach. This approach was used in the earlier works \cite{bk:MkhizeMoyoMeleshko[2014], bk:PopovychKunzingerEshraghi[2010],bk:Kasatkin[2012],%
bk:GrigorievMeleshkoSuriyawichitseranee[2013]} and references therein. The algebraic algorithm assists in simplifying the study of the group classification.

The rest of the paper is organized as follows:
Section $2$ of the paper gives a background study of systems of nonlinear equations and describes the concept of a reducible system and irreducible systems. Furthermore the equivalence transformations and determining equations are discussed.
The present paper focuses
on irreducible systems. Section $3$ discusses linear equations where the coefficient matrix $B(x)$ defined in this section is reduced to zero and the property of the trace of the coefficient matrix $C(x)$ is also discussed. The determining equations and the commutator tables are computed. Section $4$ gives a strategy and detailed approach for the group classification that uses the algebraic approach using the optimal system of subalgebras of the Lie algebras for group classification. Some interesting observations are made. The main thrust of the paper is in Section $5$ where the classes of systems admitting the associated Lie algebras obtained in Section $4$ are delineated. The relationship between the classification of systems of second-order ordinary differential equations with the classification of the coefficient matrix $A$ reducing it to the Jordan form is noted. In this Section also lies the main Theorem and results. Section $6$ is the discussion on systems admitting generators of the form $X_{A}$ defined in Section $5$. Finally Section $7$ gives the conclusion of the paper and Section $8$ is the Appendix where a detailed analysis of matrix equations is given.

\section{Background study of systems of the
form $\mathbf{y}^{\prime\prime}=\mathbf{F}(x,\mathbf{y})$}

We give a preliminary study of systems of nonlinear equations by considering a system of second-order ordinary differential equations of the form

\begin{equation}
\mathbf{y}^{\prime\prime}=\mathbf{F}(x,\mathbf{y}),\label{eq:main}
\end{equation}
 where
\[
\mathbf{y}=\left(\begin{array}{c}
y_{1}\\
y_{2}\\
...\\
y_{m}
\end{array}\right),\ \ \ \mathbf{F}=\left(\begin{array}{c}
F_{1}\\
F_{2}\\
...\\
F_{m}
\end{array}\right).
\]

\subsection{Equivalence transformations}

System (\ref{eq:main}) has the following equivalence transformations:

(a) a linear change of the dependent variables $\widetilde{\mathbf{y}}=P\mathbf{y}$
with a constant nonsingular $m\times m$ matrix;

(b) the change \[\widetilde{y}_{i}=y_{i}+\varphi_{i}(x),\quad(i=1,2,..,m);
\]

(c) a transformation related with the change
\[
\widetilde{x}=\varphi(x),\:\widetilde{y}_{i}=y_{i}\psi(x),\:(i=1,2,...,m),
\]
 where the functions $\varphi(x)$ and $\psi(x)$ satisfy the condition
\begin{equation}
\frac{\varphi^{\prime\prime}}{\varphi^{\prime}}=2\frac{\psi^{\prime}}{\psi}.\label{eq:feb27}
\end{equation}

We call a system of equations (\ref{eq:main}) reducible if it is
equivalent to a system which has a proper subsystem of fewer dimension
or if it is equivalent with respect to a change of the dependent and independent
variables to a linear system
\begin{equation}
\mathbf{y}^{\prime\prime}=C\mathbf{y},\label{eq:mainlinear-1}
\end{equation}
 where $C$ is a constant matrix. In the present paper irreducible
systems are considered.

\subsection{Determining equations}

Determining equations for irreducible systems in matrix form are given by
\begin{equation}
2\xi\mathbf{F}_{x}+3\xi^{\prime}\mathbf{F}+\left(((A+\xi^{\prime}E)\mathbf{y}+\zeta)\cdot\nabla\right)\mathbf{F}-A\mathbf{F}=\xi^{\prime\prime\prime}\mathbf{y}+\zeta^{\prime\prime},\label{eq:aug1011.1}
\end{equation}
 where the matrix $A=(a_{ij})$ is constant. The associated infinitesimal generator
is
\[
X=2\xi(x)\partial_{x}+(A\mathbf{y}+\zeta(x))\cdot\nabla,
\]
where $\nabla=\left(\partial_{y_{1}},\partial_{y_{2}},...,\partial_{y_{m}}\right)^{t}$. Here "$\cdot$" means the scalar product:
$\mathbf{b}\cdot\nabla=b_{i}\partial_{y_{i}}$ and the standard
agreement, summation with respect to a repeat index, is used.

Applying the change $\widetilde{\mathbf{y}}=P\mathbf{y}$ where $P$
is a nonsingular $m\times m$ matrix with constant entries, equations
(\ref{eq:main}) become
\[
\widetilde{\mathbf{y}}^{\prime\prime}=\widetilde{\mathbf{F}}(x,\widetilde{\mathbf{y}})
\]
 with
\[
\widetilde{\mathbf{F}}(x,\widetilde{\mathbf{y}})=P\mathbf{F}(x,P^{-1}\widetilde{\mathbf{y}}).
\]
 The partial derivatives with respect to the variables $\mathbf{y}$
are changed as follows:
\[
\mathbf{h}\cdot\nabla=(P\mathbf{h})\cdot\widetilde{\nabla}.
\]
 Hence equations (\ref{eq:aug1011.1}) become
\[
2\xi\widetilde{\mathbf{F}}{}_{x}+3\xi^{\prime}\widetilde{\mathbf{F}}
+\left(((\widetilde{A}+\xi^{\prime}E)\widetilde{\mathbf{y}}
+\widetilde{\zeta})\cdot\widetilde{\nabla}\right)\widetilde{\mathbf{F}}
-\widetilde{A}\widetilde{\mathbf{F}}-\xi^{\prime\prime\prime}\widetilde{\mathbf{y}}
+\widetilde{\zeta}^{\prime\prime}=0,
\]
 where
\[
\widetilde{A}=PAP^{-1},\,\,\,\widetilde{\zeta}=P\zeta.
\]
 This means that the equivalence transformation $\widetilde{\mathbf{y}}=P\mathbf{y}$
reduces equation (\ref{eq:aug1011.1}) to the same form with the matrix
$A$ and the vector $\zeta$ changed. The infinitesimal generator
is also changed as follows:
\[
X=2\xi\partial_{x}+(\widetilde{A}\widetilde{\mathbf{y}}+\widetilde{\zeta})\cdot\widetilde{\nabla}.
\]

\section{Systems of linear equations}

Systems of linear second-order ordinary differential equations have
the following form,
\begin{equation}
\mathbf{y}^{\prime\prime}=B(x)\mathbf{y}^{\prime}+C(x)\mathbf{y}+\mathbf{f}(x)\label{eq:mainlinear}
\end{equation}
 where $B(x)$ and $C(x)$ are matrices, and $f(x)$ is a vector.
Using a particular solution $\mathbf{y}_{p}(x)$ and the change
\[
\mathbf{y}=\widetilde{\mathbf{y}}+\mathbf{y}_{p},
\]
we can assume that $\mathbf{f}(x)=0$ without loss of generality.
The matrix $B(x)$ or $C(x)$ can also be assumed to be zero if the
change, $\mathbf{y}=H(x)\widetilde{\mathbf{y}}$, where $H=H(x)$ is a nonsingular matrix, is used.
In the current paper the matrix $B(x)$ is
reduced to zero. In this case the function ${\bf F}$ in equation (\ref{eq:main})
is a linear function of $\mathbf{y}$:
\[
\mathbf{F}(x,\mathbf{y})=C(x)\mathbf{y}.
\]

Any linear system of second-order ordinary differential equations
\begin{equation}
\mathbf{y}^{\prime\prime}=C(x)\mathbf{y}\label{tg1}
\end{equation}
 admits the set of trivial generators
\[
\mathbf{y}\cdot\nabla,\,\,\,\mathbf{h}(x)\cdot\nabla,
\]
 where $\mathbf{h}^{\prime\prime}=C\mathbf{h}$.

Excluding the trivial generators, the determining equations (\ref{eq:aug1011.1}),
after their splitting with respect to $\mathbf{y}$, become
\begin{equation}
2\xi C^{\prime}+CA-AC=\xi^{\prime\prime\prime}E-4\xi^{\prime}C,\label{eq:oct5.2013.20}
\end{equation}
 where $E$ is the unit $m\times m$ matrix, and the admitted generator
has the form
\[
X=2\xi\partial_{x}+((A+\xi^{\prime}E)\mathbf{y})\cdot\nabla.
\]

\subsection{Simplifications of systems}

Applying the change of the dependent and independent variables
\begin{equation}
\tilde{x}=\varphi(x),\ \ \ \mathbf{\tilde{y}}=\psi(x)\mathbf{y}\label{eq:aug14.1}
\end{equation}
 satisfying the condition
\begin{equation}
\frac{\varphi^{\prime\prime}}{\varphi^{\prime}}=2\frac{\psi^{\prime}}{\psi},\label{eq:apr28.3}
\end{equation}
 system (\ref{tg1}) becomes
\begin{equation}
\mathbf{\tilde{y}}^{\prime\prime}=\tilde{C}\mathbf{\tilde{y}},\label{eq:apr27.1}
\end{equation}
 where
\[
\tilde{C}=\varphi^{\prime}{}^{-2}\left(C-\frac{\rho^{\prime\prime}}{\rho}E\right),\ \ \rho=\frac{1}{\psi}.
\]
The group classification problem usually becomes simpler after reducing
the number of arbitrary elements. In order to reduce the number of entries
of the matrix $\tilde{C}$ one can choose the function $\psi$ such
that%
\footnote{ This change was used in \cite{bk:WafoMahomed[2000]} for the case
of $m=2$%
} $tr(\tilde{C})=0$. This condition leads to the equation
\begin{equation}
\rho^{\prime\prime}-\frac{tr(C)}{m}\rho=0.\label{eq:apr28.2}
\end{equation}
 Notice that a transformation of the form (\ref{eq:aug14.1}) with
\begin{equation}
\psi=\alpha(x+\beta)^{-1}\label{eq:dec02.13.1}
\end{equation}
conserves the property $tr(C)=0$. Here $\alpha$ and $\beta$ are
constants. In fact, if $tr(C)=0$, then because of $\rho^{\prime\prime}=0$
we have that $tr(\tilde{C})=0$. In particular, the equivalence transformation
with
\begin{equation}
\varphi=\psi=\frac{1}{x}\label{eq:dec19.13.1}
\end{equation}
is an involution.

\subsection{Determining equations}

As noted earlier, for the group classification we can assume that $tr(C)=0$.
Taking the trace in (\ref{eq:oct5.2013.20}), one finds that $\xi^{\prime\prime\prime}=0$
or
\[
\xi=\frac{1}{2}(k_{1}x^{2}+k_{3})+k_{2}x,
\]
 where $k_{i}$, ($i=0,1,2$) are constants. Hence, nontrivial admitted generators
take the form
\[
X=k_{1}X_{1}+k_{2}X_{2}+k_{3}X_{3}+X_{A},
\]
 where
\[
X_{1}=x(x\partial_{x}+\mathbf{y}\cdot\nabla),\,\,\, X_{2}=2x\partial_{x}+\mathbf{y}\cdot\nabla,\,\,\, X_{3}=\partial_{x},\,\,\, X_{A}=(A\mathbf{y})\cdot\nabla.
\]
Notice that the generator $X_{2}$ can be simplified by subtracting the
trivial admitted generator $\mathbf{y}\cdot\nabla$. However, we keep
it in the presented form due to the simplicity of the commutator
\[
[X_{1},X_{3}]=-X_{2}.
\]

 The determining equations become
\begin{equation}
(k_{1}x^{2}+2k_{2}x+k_{3})C^{\prime}+CA-AC+4(k_{1}x+k_{2})C=0.\label{eq:oct5.2013.22}
\end{equation}
Thus we find that an admitted Lie algebra of nontrivial generators
is composed by the generators $X_{1}$, $X_{2}$, $X_{3}$ and $X_{A}$.

To study the problem further we need to construct the commutator table of the generators
$X_{1}$, $X_{2}$, $X_{3}$:
\[
\begin{array}{c|ccc}
\hline  & X_{1} & X_{2} & X_{3}\\
\hline X_{1} & 0 & -2X_{1} & -X_{2}\\
X_{2} & 2X_{1} & 0 & -2X_{3}\\
X_{3} & X_{2} & 2X_{3} & 0
\end{array}
\]

\section{Strategy for the group classification}

One of the methods for analyzing relations between the constants and
undefined functions consists of employing the algorithm developed
for the gas dynamics equations \cite{bk:Ovsiannikov[1978]}. This
algorithm allows one to study all possible admitted Lie algebras without
omission. Unfortunately, it is difficult to implement for system (\ref{tg1}).
Observe also that sometimes in this approach it is difficult to select out equivalent
cases with respect to equivalence transformations.

In \cite{bk:MkhizeMoyoMeleshko[2014],bk:PopovychKunzingerEshraghi[2010],bk:Kasatkin[2012],%
bk:GrigorievMeleshkoSuriyawichitseranee[2013]}%
\footnote{See also references therein.}
a different approach was applied for group classification. We call
this approach an algebraic algorithm. In most applications the algebraic
algorithm essentially reduces the study of group classification to
a simpler problem. The reduction occurs because the process of solving
determining equations is split into two steps, where on the first
step the constants of admitted generators are defined using the property
for admitted generators to compose a Lie algebra. In the present paper we follow the algebraic approach.

\subsection{Relations between automorphisms and equivalence transformations}

Generators of admitted Lie algebras have the form
\[
X=x_{1}X_{1}+x_{2}X_{2}+x_{3}X_{3}+X_{A}.
\]
 The commutator of two generators
\[
X=x_{1}X_{1}+x_{2}X_{2}+x_{3}X_{3}+X_{A_{1}},\,\,\, Z=z_{1}X_{1}+z_{2}X_{2}+z_{3}X_{3}+X_{A_{2}}
\]
 is
\[
[X,Z]=\alpha X_{1}+\beta X_{2}+\gamma X_{3}+((A_{2}A_{1}-A_{1}A_{2})\mathbf{y})\cdot\nabla,
\]
 where
\[
\alpha=-2(x_{1}z_{2}-x_{2}z_{1}),\ \ \ \beta=-(x_{1}z_{3}-x_{3}z_{1}),\ \ \ \gamma=-2(x_{2}z_{3}-x_{3}z_{2}).
\]
 Hence, one can notice that the first part of the admitted generators
$x_{1}X_{1}+x_{2}X_{2}+x_{3}X_{3}$ is a subalgebra of the Lie
algebra\footnote{This Lie algebra corresponds to the
algebra type YIII in the Bianchi classification.}
$L_{3}=\{X_{1},X_{2},X_{3}\}$. Recall that all nonequivalent subalgebras
with respect to automorphisms present an optimal system
of subalgebras.
An optimal system of subalgebras of the algebra type VIII in the
Bianchi classification was performed in \cite{bk:PateraWinternitz[1977]}.

We further show that the action of equivalence transformations
conserving the property $tr(C)=0$ is similar to the action of automorphisms.
This property allows one to use an optimal system of subalgebras of
the Lie algebra $L_{3}$ for group classification.

In fact, automorphisms of $L_{3}$ are
\[
Aut_{1}:\,2x_{2}\partial_{x_{1}}+x_{3}\partial_{x_{2}}\,\,\,\,\,\,\,\,\,\,\bar{x}_{1}=x_{1}+2ax_{2}+a^{2}x_{3},\,\,\bar{x}_{2}=x_{2}+ax_{3};
\]
\[
Aut_{2}:\,\,\,\,\,\,\,\,\,\,\,\bar{x}_{1}=x_{1}e^{a},\,\,\bar{x}_{3}=x_{3}e^{-a};
\]
\[
Aut_{3}:\, x_{1}\partial_{x_{2}}+2x_{2}\partial_{x_{3}}\,\,\,\,\,\,\,\,\,\,\,\,\bar{x}_{2}=x_{2}+ax_{1},\,\,\bar{x}_{3}=x_{3}+2ax_{2}+a^{2}x_{1}.
\]
 Here and further on only changeable coordinates of the generator
are presented. The equivalence transformation (\ref{eq:aug14.1})
with
\[
\varphi=\frac{x}{1-ax},\,\,\,\psi=(x+a)^{-1}
\]
 changes the coordinates as follows
\[
(x_{1}X_{1}+x_{2}X_{2}+x_{3}X_{3})(\varphi(x))\partial_{\bar{x}}=\bar{x}^{2}(x_{3}a^{2}+2ax_{2}+x_{1})\partial_{\bar{x}}+2\bar{x}(x_{2}+ax_{3})\partial_{\bar{x}}+x_{3}\partial_{\bar{x}}
\]
 or
\[
x_{1}X_{1}+x_{2}X_{2}+x_{3}X_{3}=(x_{3}a^{2}+2ax_{2}+x_{1})\bar{X}_{1}+(x_{2}+ax_{3})\bar{X}_{2}+x_{3}\bar{X}_{3}.
\]
 Hence, this equivalence transformation is similar to the automorphism
$Aut_{1}$. The equivalence transformation $\bar{x}=ax$ is equivalent
to the automorphism $Aut_{2}$. The equivalence transformation $\bar{x}=x-a$
corresponding to the shift of $x$ is similar to the automorphism
$Aut_{3}:$
\[
x_{1}X_{1}+x_{2}X_{2}+x_{3}X_{3}=x_{1}\bar{X}_{1}+(x_{2}+ax_{1})\bar{X}_{2}+(x_{3}+a^{2}x_{1}+2ax_{2})\bar{X}_{3}.
\]

The use of the optimal system of subalgebras of $L_{3}$ for the group
classification is similar to the two-step algorithm of constructing
an optimal system of subalgebras \cite{bk:Ovsiannikov[1993opt]},
where on the first step an optimal system of a fewer dimension subalgebras
is constructed.


As mentioned above an optimal system of subalgebras of the Lie algebra $L_{3}$ was studied in \cite{bk:PateraWinternitz[1977]} and it consists of the list:
\begin{equation}
\label{eq:dec2613.1_0}
\begin{array}{cl}
1. & X_{2};\\
2. & X_{3};\\
3. & X_{1}+X_{3};\\
4. & X_{2},\, X_{3};\\
5. & X_{1},\, X_{2},\, X_{3}.
\end{array}
\end{equation}

\subsection{Classes of systems admitting generators with $\xi\neq0$}
\label{sec_4.3}

Using the optimal system of subalgebras (\ref{eq:dec2613.1_0}), we
can conclude that all systems of linear second-order ordinary differential
equations (\ref{tg1}) admitting generators with $\xi\neq0$ are separated
into the classes admitting the following Lie algebras:
\begin{equation}
\begin{array}{cl}
1. & X_{2}+X_{A_{2}};\\
2. & X_{3}+X_{A_{3}};\\
3. & X_{1}+ X_{3}+X_{A_{1}};\\
4. & X_{2}+X_{A_{2}},\, X_{3}+X_{A_{3}};\\
5. & X_{1}+X_{A_{1}},\, X_{2}+X_{A_{2}},\, X_{3}+X_{A_{3}}.
\end{array}\label{eq:dec2613.1}
\end{equation}
 Here the numeration of matrices $A_{i}$ is used to ease the tracking of
their relations with the generators $X_{i}$.

Except the generators presented in (\ref{eq:dec2613.1}), systems
(\ref{tg1}) can admit several generators of the form $X_{A}$. The
determining equations in this case are
\[
CA-AC=0.
\]
 According to the algebraic study considered in the Appendix, the number
of equations has to be even ($m=2n$), and excluding the trivial generator
$\mathbf{y}\cdot\nabla$, the matrix $A$ for one of the admitted
generators of the form $X_{A}$, can be chosen as $A=B_{d_{1}}$, where
\[
B_{d_{1}}=\left(\begin{array}{cccc}
B_{1} & 0 & ... & 0\\
0 & B_{1} & ... & 0\\
... & ... & ... & ...\\
0 & 0 & ... & B_{1}
\end{array}\right),\ \ \ B_{1}=\left(\begin{array}{cc}
0 & 1\\
-1 & 0
\end{array}\right).
\]
 The latter matrices have the properties
\[
B_{1}^{-1}=-B_{1},\,\,\, B_{d_{1}}^{-1}=-B_{d_{1}}.
\]
Notice
also that in this case the matrix $C$ consists of blocks of the form
\[
C_{ij}=\left(\begin{array}{cc}
\alpha_{ij} & \beta_{ij}\\
-\beta_{ij} & \alpha_{ij}
\end{array}\right),
\]
 and the determining equations (\ref{eq:oct5.2013.22}) for this generator
$(B_{d_{1}}\mathbf{y})\cdot\nabla$ are reduced to the equations
\begin{equation}
(k_{1}x^{2}+2k_{2}x+k_{3})\alpha_{ij}^{\prime}+4(k_{1}x+k_{2})\alpha_{ij}=0,\,\,\,(k_{1}x^{2}+2k_{2}x+k_{3})\beta_{ij}^{\prime}+4(k_{1}x+k_{2})\mathbf{\beta}_{ij}=0.\label{eq:oct5.2013.35}
\end{equation}

In addition we also conclude here that for systems with an odd number of dependent
variables and having  nontrivial admitted generators, only these
admitted Lie algebras (\ref{eq:dec2613.1}) are possible.

\section{Classes of systems admitting Lie algebras (\ref{eq:dec2613.1})}

In this section we consider systems corresponding to the Lie algebras presented
in (\ref{eq:dec2613.1}).

\subsection{Systems admitting the generator $X_{2}+X_{A}$ }

For simplifying the determining equations in this case we apply
the change (\ref{eq:aug14.1}) with
\[
\varphi=\ln(x),\,\,\psi=x^{-1/2}.
\]
 The determining equations
\[
2xC^{\prime}+CA-AC+4C=0,
\]
 become
\[
2\frac{d}{d\tilde{x}}\tilde{C}+\tilde{C}A-A\tilde{C}=0,\,\,\, tr(\tilde{C})=\frac{m}{4}.
\]
 Thus,
\[
\tilde{C}=e^{\tilde{x}\tilde{A}}C_{0}e^{-\tilde{x}\tilde{A}},
\]
 where $\tilde{A}=\frac{1}{2}A,$ and $C_{0}$ is an arbitrary matrix
with $tr(C_{0})=m/4$. It is also assumed that
\[
AC_{0}-C_{0}A\neq0,
\]
 because otherwise $AC-CA=0$ which implies that the matrix $C$
is constant. Note that because $\frac{d}{d\tilde{x}}(tr(\tilde{C}))=0$,
we have that $tr(\tilde{C})=m/4$. The admitted generator is
\[
X=\partial_{\tilde{x}}+(\tilde{A}\mathbf{\tilde{y}})\cdot\tilde{\nabla}.
\]
Here the part related with the trivial admitted generator $\mathbf{\tilde{y}}\cdot\tilde{\nabla}$
is omitted.

\textbf{Remark}. Further classification of systems of second-order
ordinary differential equations of this type is related with the classification
of the matrix $A$, reducing it to one of Jordan forms.
This remark applies to other cases discussed further on.

\subsection{Systems admitting the generator $X_{3}+X_{A}$ }

In this case the general solution of the determining equations is
\[
C=e^{-xA}C_{0}e^{-xA},
\]
 where $C_{0}$ is an arbitrary matrix with $tr(C_{0})=0$. It is observed that
because $\frac{d}{d\tilde{x}}(tr(C))=0$, we have that $tr(C)=0$.
The admitted generator is
\[
X=\partial_{x}+(A\mathbf{y})\cdot\nabla.
\]

\subsection{Systems admitting the generator $X_{1}+ X_{3}+X_{A}$ }

For the generator
\[
X_{1}+ X_{3}+X_{A}=(x^{2}+1)\partial_{x}+(A\mathbf{y})\cdot\nabla
\]
 the determining equations are
\[
(x^{2}+1)C^{\prime}+CA-AC+4xC=0.
\]
To simplify the determining equations we apply the change
(\ref{eq:aug14.1}) with
\[
\varphi^{\prime}=(x^{2}+1)^{-1},\,\,\psi=(x^{2}+1)^{-1/2}.
\]
 The determining equations become
\[
\frac{d}{d\tilde{x}}\tilde{C}+\tilde{C}A-A\tilde{C}=0,\,\,\, tr(\tilde{C})=-1.
\]
 Hence,
\[
\tilde{C}=e^{-\tilde{x}A}C_{0}e^{\tilde{x}A},
\]
 where $C_{0}$ is an arbitrary matrix with $tr(C_{0})=-1$.
We point out that because $\frac{d}{d\tilde{x}}(tr(\tilde{C}))=0$, one
also has that $tr(\tilde{C})=-1$. The admitted generator is
\[
X=\partial_{\tilde{x}}+(A\mathbf{\tilde{y}})\cdot\tilde{\nabla}.
\]

\subsection{Discussion on systems admitting one-dimensional Lie algebras from
(\ref{eq:dec2613.1}) }

The study above allows us to conclude that irreducible linear systems
(\ref{tg1}) admitting a generator with $\xi\neq0$ are equivalent
to a system (\ref{tg1}) where
\[
C(x)=e^{xA}C_{0}e^{-xA},
\]
 and $C_{0}A-AC_{0}\neq0$. The admitted generator is
\[
X=\partial_{x}+(A\mathbf{y})\cdot\nabla.
\]
There is no necessity to take care on $tr(C_{0})$. It was only necessary
for being sure that none of the linear systems admitting a Lie group is
missed.

\subsection{Systems admitting the generators $X_{2}+X_{A_{2}}$ and \textmd{\normalsize $X_{3}+X_{A_{3}}$}}

The commutator of these generators is
\[
[X_{2}+(A_{2}\mathbf{y})\cdot\nabla,X_{3}+(A_{3}\mathbf{y})\cdot\nabla]=-2(X_{3}+(A_{3}\mathbf{y})\cdot\nabla)+((A_{3}A_{2}-A_{2}A_{3}+2A_{3})\mathbf{y})\cdot\nabla.
\]
 Since the admitted Lie algebra is two-dimensional, then
\begin{equation}
A_{3}(A_{2}+2E)-A_{2}A_{3}=0.\label{eq:dec8.13.5-2}
\end{equation}

As noted above, because the generator $X_{3}+X_{A_{3}}$ is admitted,
then
\[
C^{\prime}=A_{3}C-CA_{3},
\]
 and then
\[
C=e^{xA_{3}}C_{0}e^{-xA_{3}},
\]
 where $C_{0}$ is an arbitrary matrix with $tr(C_{0})=0$, and it
is also assumed that
\[
A_{3}C_{0}-C_{0}A_{3}\neq0.
\]
 The determining equations for the generator $X_{2}+X_{A_{2}}$ are
\begin{equation}
-2x(CA_{3}-A_{3}C)+CA_{2}-A_{2}C+4C=0.\label{eq:dec813.1-3}
\end{equation}
Note that the substitution in this equation $x=0$ gives
\[
A_{2}C_{0}-C_{0}(A_{2}+4E)=0.
\]
 Multiplying from the left hand side by $e^{xA_{3}}$ and by $e^{-xA_{3}}$
from the right hand side, equations (\ref{eq:dec813.1-3}) can be rewritten as
\begin{equation}
2x(A_{3}C_{0}-C_{0}A_{3})+4C_{0}+C_{0}e^{-xA_{3}}A_{2}e^{xA_{3}}-e^{-xA_{3}}A_{2}e^{xA_{3}}C_{0}=0.\label{eq:dec613.1-3}
\end{equation}
 Differentiating (\ref{eq:dec613.1-3}) and using the property that
the matrices $e^{xA_{3}}$ and $e^{-xA_{3}}$ commute with the
matrix $A_{3}$, we obtain
\begin{equation}
2(A_{3}C_{0}-C_{0}A_{3})+C_{0}e^{-xA_{3}}(-A_{3}A_{2}+A_{2}A_{3})e^{xA_{3}}-e^{-xA_{3}}(-A_{3}A_{2}+A_{2}A_{3})e^{xA_{3}}C_{0}=0.\label{eq:dec613.2-3}
\end{equation}
Equation (\ref{eq:dec613.2-3}) can be rewritten as
\[
CB-BC=0,
\]
 where
\begin{equation}
B=A_{2}A_{3}-A_{3}(A_{2}+2E).\label{eq:dec6.13.6-3}
\end{equation}
Due to (\ref{eq:dec8.13.5-2}) we have that $B=0$. Hence we obtain a system of algebraic equations for the matrices $C_{0}$,
$A_{2}$ and $A_{3}$ given as follows:
\begin{equation}
A_{2}C_{0}-C_{0}(A_{2}+4E)=0,
\label{eq:dec6.13.10-3}
\end{equation}
\begin{equation}
A_{2}A_{3}-A_{3}(A_{2}+2E)=0,
\label{eq:dec6.13.11-3}
\end{equation}
\begin{equation}
tr(C_{0})=0,\,\,\, A_{3}C_{0}-C_{0}A_{3}\neq0.
\label{eq:dec6.13.12-3}
\end{equation}

Since $C_{0}\neq0$, then equation (\ref{eq:dec6.13.10-3}) gives us
that the set of eigenvalues of the matrix $A_{2}$ has intersections
with the set of eigenvalues of the matrix%
\footnote{Discussion of solving matrix equations can be found in \cite{bk:Gantmacher}.%
} $A_{2}+4E$. Since $A_{3}\neq\beta E$ (for any $\beta$), then equation
(\ref{eq:dec6.13.11-3}) gives us that the set of eigenvalues of the
matrix $A_{2}$ has intersections with the set of eigenvalues of the
matrix $A_{2}+2E$. Hence there exists a number $\lambda$ such that
$\lambda$, $\lambda+2$ and $\lambda+4$ are eigenvalues of the matrix
$A_{2}$.

For example, if $m=3$, one can assume that
\[
A_{2}=\left(\begin{array}{ccc}
0 & 0 & 0\\
0 & 2 & 0\\
0 & 0 & 4
\end{array}\right).
\]
 Since the matrix $A_{2}$ is diagonal, the general solution of equation
(\ref{eq:dec6.13.10-3}) and (\ref{eq:dec6.13.11-3}) is trivially
obtained (see in \cite{bk:Gantmacher}):
\[
C_{0}=\left(\begin{array}{ccc}
0 & 0 & 0\\
0 & 0 & 0\\
c & 0 & 0
\end{array}\right),\,\,\, A_{3}=\left(\begin{array}{ccc}
0 & 0 & 0\\
a & 0 & 0\\
0 & b & 0
\end{array}\right),
\]
 where $a$, $b$ and $c$ are constant. Since conditions (\ref{eq:dec6.13.12-3})
are not satisfied for these matrices, then there are no such generators
in the case of $m=3$. This is also supported by the study in \cite{bk:SuksernMoyoMeleshko[2014]}.

Let us also consider as example the case where $m=4$. One can assume
that the matrix $A_{2}$ is one of the matrices:
\[
\left(\begin{array}{cccc}
a & 0 & 0 & 0\\
0 & 0 & 0 & 0\\
0 & 0 & 2 & 0\\
0 & 0 & 0 & 4
\end{array}\right),\,\,\,\left(\begin{array}{cccc}
0 & 1 & 0 & 0\\
0 & 0 & 0 & 0\\
0 & 0 & 2 & 0\\
0 & 0 & 0 & 4
\end{array}\right),\,\,\,\left(\begin{array}{cccc}
0 & 0 & 0 & 0\\
0 & 2 & 1 & 0\\
0 & 0 & 2 & 0\\
0 & 0 & 0 & 4
\end{array}\right),\,\,\,\left(\begin{array}{cccc}
0 & 0 & 0 & 0\\
0 & 2 & 0 & 0\\
0 & 0 & 4 & 1\\
0 & 0 & 0 & 4
\end{array}\right).
\]
 Calculations show that all of these cases are reduced to a reducible
system: either there is a subsystem with fewer dimension or the matrix
$C$ is constant. Thus for $m=4$ there is no such admitted Lie subalgebra.

\subsection{Systems admitting the generators $X_{1}+X_{A_{1}}$, $X_{2}+X_{A_{2}}$
and \textmd{\normalsize $X_{3}+X_{A_{3}}$}}

The commutators of these generators are
\[
\begin{array}{rcl}
[X_{1}+X_{A_{1}},X_{2}+X_{A_{2}}] & = & -2(X_{1}+X_{A_{1}})+((A_{2}A_{1}-A_{1}A_{2}+2A_{1})\mathbf{y})\cdot\nabla,\\
\ [X_{1}+X_{A_{1}},X_{3}+X_{A_{3}}] & = & -(X_{2}+X_{A_{2}})+((A_{3}A_{1}-A_{1}A_{3}+A_{2})\mathbf{y})\cdot\nabla,\\
\ [X_{2}+X_{A_{2}},X_{3}+X_{A_{3}}] & = & -2(X_{3}+X_{A_{3}})+((A_{3}A_{2}-A_{2}A_{3}+2A_{3})\mathbf{y})\cdot\nabla.
\end{array}
\]
 Hence,
\begin{equation}
A_{1}A_{2}-A_{2}A_{1}=2A_{1},\,\,\, A_{1}A_{3}-A_{3}A_{1}=A_{2},\label{eq:jan14.10.3}
\end{equation}
\begin{equation}
A_{2}A_{3}-A_{3}A_{2}=2A_{3}.\label{eq:jan14.10.4}
\end{equation}
 The determining equations are
\begin{equation}
C^{\prime}+CA_{3}-A_{3}C=0,\ \ \ 2xC^{\prime}+CA_{2}-A_{2}C+4C=0,\label{eq:jan14.10.1}
\end{equation}
\begin{equation}
x^{2}C^{\prime}+CA_{1}-A_{1}C+4xC=0.\label{eq:jan14.10.2}
\end{equation}

As shown in the previous section, the general solution of equations
(\ref{eq:jan14.10.1}), (\ref{eq:jan14.10.2}) and (\ref{eq:jan14.10.4})
is
\[
C=e^{xA_{3}}C_{0}e^{-xA_{3}},
\]
 where the matrices $C_{0}$, $A_{2}$ and $A_{3}$ satisfy conditions
(\ref{eq:dec6.13.10-3}), (\ref{eq:dec6.13.11-3}) and (\ref{eq:dec6.13.12-3}).

The remaining determining equations (\ref{eq:jan14.10.2}) become
\[
S_{1}=0,
\]
 where
\[
S_{1}=4xC_{0}-x^{2}(C_{0}A_{3}-A_{3}C_{0})+C_{0}e^{-xA_{3}}A_{1}e^{xA_{3}}-e^{-xA_{3}}A_{1}e^{xA_{3}}C_{0}.
\]
 Note that
\[
S_{1}^{\prime}=4C_{0}-2x(C_{0}A_{3}-A_{3}C_{0})+C_{0}e^{-xA_{3}}(A_{1}A_{3}-A_{3}A_{1})e^{xA_{3}}-e^{-xA_{3}}(A_{1}A_{3}-A_{3}A_{1})e^{xA_{3}}C_{0}
\]
 and
\[
\begin{array}{rl}
e^{xA_{3}}S_{1}^{\prime\prime}e^{-xA_{3}}= & -2(CA_{3}-A_{3}C)+C[(A_{1}A_{3}-A_{3}A_{1})A_{3}-A_{3}(A_{1}A_{3}-A_{3}A_{1})]\\
 & -[(A_{1}A_{3}-A_{3}A_{1})A_{3}-A_{3}(A_{1}A_{3}-A_{3}A_{1})]C.
\end{array}
\]
Conditions (\ref{eq:dec6.13.10-3}), (\ref{eq:jan14.10.3}) and (\ref{eq:jan14.10.4})
imply  that
\[
S_{1}^{\prime}(0)=0,\ \ \ S_{1}^{\prime\prime}=0.
\]
 Hence
\[
S_{1}(x)=S_{1}(0)=C_{0}A_{1}-A_{1}C_{0}=0.
\]
 Thus, we obtain the following conditions for the matrices $C_{0}$,
$A_{1}$, $A_{2}$ and $A_{3}$:
\begin{equation}
A_{1}A_{2}-A_{2}A_{1}=2A_{1},\,\,\, A_{1}A_{3}-A_{3}A_{1}=A_{2},\,\,\, A_{2}A_{3}-A_{3}A_{2}=2A_{3},\label{eq:jan14.11.10}
\end{equation}
\begin{equation}
A_{2}C_{0}-C_{0}(A_{2}+4E)=0,\,\,\, C_{0}A_{1}-A_{1}C_{0}=0.\label{eq:jan14.11.11}
\end{equation}

{\bf Remark}. Defining from these equations $A_{2}=A_{1}A_{3}-A_{3}A_{1}$, and
substituting it into the remaining equations (\ref{eq:jan14.11.10})
and (\ref{eq:jan14.11.11}) we obtain only equations for the matrices
$C_{0}$, $A_{1}$ and $A_{3}$:
\begin{equation}
A_{1}^{2}A_{3}-2A_{1}A_{3}A_{1}+A_{3}A_{1}^{2}=2A_{1},\,\,\, A_{1}A_{3}^{2}-2A_{3}A_{1}A_{3}+A_{3}^{2}A_{1}=2A_{3},\label{eq:jan14.11.12}
\end{equation}
\begin{equation}
(A_{1}A_{3}-A_{3}A_{1})C_{0}-C_{0}(A_{1}A_{3}-A_{3}A_{1})-4C_{0}=0,\,\,\, C_{0}A_{1}-A_{1}C_{0}=0.\label{eq:jan14.11.13}
\end{equation}

\subsection{Summary of the results}

We note that if  system (\ref{tg1}) admits one generator with
$\xi\neq0$, then without loss of generality one can assume that $\xi=1$.

As a result of this section we derive the Theorem.

\textbf{Theorem}. Irreducible linear systems (\ref{tg1}) admitting
one- two- or three-dimensional Lie algebras (\ref{eq:dec2613.1})
are equivalent to one of the following cases.

(a) For one-dimensional Lie algebras
\[
C(x)=e^{xA}C_{0}e^{-xA},
\]
 and $C_{0}A-AC_{0}\neq0$. The admitted generator is
\[
X=X_{3}+(A\mathbf{y})\cdot\nabla.
\]

(b) For two-dimensional Lie algebras
\[
X_{2}+X_{A_{2}},\,\,\, X_{3}+X_{A_{3}}
\]
 the system (\ref{tg1}) has
\[
C(x)=e^{xA_{3}}C_{0}e^{-xA_{3}},
\]
 where the matrices $C_{0}$, $A_{2}$ and $A_{3}$ satisfy the conditions:
\[
A_{2}A_{3}-A_{3}(A_{2}+2E)=0,
\]
 and
\[
A_{2}C_{0}-C_{0}(A_{2}+4E)=0,\,\,\, tr(C_{0})=0,\,\,\, A_{3}C_{0}-C_{0}A_{3}\neq0.
\]

(c) For three-dimensional Lie algebras
\[
X_{1}+X_{A_{1}},\,\,\, X_{2}+X_{A_{2}},\,\,\, X_{3}+X_{A_{3}}
\]
 the system (\ref{tg1}) has
\[
C(x)=e^{xA_{3}}C_{0}e^{-xA_{3}},
\]
 where the matrices $C_{0}$, $A_{2}$ and $A_{3}$ satisfy the conditions:
\[
A_{1}A_{2}-A_{2}A_{1}=2A_{1},\,\,\, A_{1}A_{3}-A_{3}A_{1}=A_{2},\,\,\, A_{2}A_{3}-A_{3}A_{2}=2A_{3},
\]
\[
A_{2}C_{0}-C_{0}(A_{2}+4E)=0,\,\,\, C_{0}A_{1}-A_{1}C_{0}=0,\,\,\, tr(C_{0})=0,\,\,\, A_{3}C_{0}-C_{0}A_{3}\neq0.
\]

In particular, systems with an odd number of the dependent variables
($m=2n+1$) and having  nontrivial admitted generators are equivalent
to one of the cases presented in the Theorem.

\section{Discussion on systems admitting generators of the form $X_{A}$}

As noted earlier, for systems admitting nontrivial generators of
the form $X_{A}$ the number of the dependent variables is even ($m=2n$).
If these systems also admit a generator with $\xi\neq0$, then the
Lie algebras (\ref{eq:dec2613.1}) compose subalgebras of the admitted
Lie algebras. In this section systems admitting a single generator
of the form $X_{A}$ are considered. Since systems admitting the only
generator $X_{A}$ were considered in 
section \ref{sec_4.3}, we
study two, three and four-dimensional Lie algebras.

We mentioned in the previous section that if one of the admitted generators has
nonzero coefficients related with $\partial_{x}$, then one can assume
that this generator is $X_{3}+X_{A_{3}}$, and
\[
C(x)=e^{xA_{3}}C_{0}e^{-xA_{3}},
\]
 where $A_{3}C_{0}-C_{0}A_{3}\neq0$. Note that the determining
equations $CA-AC=0$ for the admitted generator $X_{A}$ leads to the conditions
that the equations
\begin{equation}
(k_{1}x^{2}+2k_{2}x+k_{3})(A_{3}C-CA_{3})+4(k_{1}x+k_{2})C=0\label{eq:dec27.13.10}
\end{equation}
 have the trivial solution with respect to the constants $k_{1}$,
$k_{2}$ and $k_{3}$: $k_{1}=0$, $k_{2}=0$ and $k_{3}=0$. In fact,
substituting $C$ into (\ref{eq:dec27.13.10}), and using commutativity
of $e^{xA_{3}}$, $e^{-xA_{3}}$ and $A_{3}$, we have
\[
(k_{1}x^{2}+2k_{2}x+k_{3})(A_{3}C_{0}-C_{0}A_{3})+4(k_{1}x+k_{2})C_{0}=0.
\]
 Splitting this equation with respect to $x$ leads to
\[
k_{1}(A_{3}C_{0}-C_{0}A_{3})=0,\,\,\, k_{2}(A_{3}C_{0}-C_{0}A_{3})+2k_{1}C_{0}=0,\,\,\, k_{3}(A_{3}C_{0}-C_{0}A_{3})+4k_{2}C_{0}=0.
\]
 Since $A_{3}C_{0}-C_{0}A_{3}\neq0$, we find sequentially that
$k_{1}=0$, $k_{2}=0$ and $k_{3}=0$.

\subsection{Two-dimensional admitted Lie algebra}

A basis of two-dimensional admitted Lie algebras can be chosen in the
forms
\[
X_{3}+(A_{3}\mathbf{y})\cdot\nabla,\,\,\,(A\mathbf{y})\cdot\nabla.
\]
 The commutator of these generators is
\[
[X_{3}+(A_{3}\mathbf{y})\cdot\nabla,(A\mathbf{y})\cdot\nabla]=((AA_{3}-A_{3}A)\mathbf{y})\cdot\nabla.
\]

For a two-dimensional admitted Lie algebra we find that
\begin{equation}
A_{3}A-AA_{3}=cA,\label{eq:dec913.1}
\end{equation}
 where $c$ is constant. Since $CA-AC=0$, the prohibition on the
reduction to fewer dimension allows us to choose $A=B_{d_{1}}$.

Equation (\ref{eq:dec913.1}) becomes
\begin{equation}
A_{3}B_{d_{1}}-B_{d_{1}}A_{3}=cB_{d_{1}}.\label{eq:dec10.13.1}
\end{equation}
 For analyzing the latter equation we represent the matrix $A_{3}$
in the form
\[
A_{3}=(A_{ij})=\left(\begin{array}{cccc}
A_{11} & A_{12} & ... & A_{1n}\\
A_{21} & A_{22} & ... & A_{2n}\\
... & ... & ... & ...\\
A_{n1} & A_{n2} & ... & A_{nn}
\end{array}\right),
\]
 where
\[
A_{ij}=\left(\begin{array}{cc}
a_{11}^{ij} & a_{12}^{ij}\\
a_{21}^{ij} & a_{22}^{ij}
\end{array}\right).
\]
 Multiplying (\ref{eq:dec10.13.1}) by $B_{d_{1}}$ from the right
hand side, and using the properties of the matrix $B_{d_{1}}$, one
finds that
\[
B_{d_{1}}A_{3}B_{d_{1}}^{-1}-A_{3}+cE=0.
\]
 The matrix $B_{d_{1}}A_{3}B_{d_{1}}^{-1}$ also has the block structure
\[
B_{d_{1}}A_{3}B_{d_{1}}^{-1}=(B_{1}A_{ij}B_{1}^{-1}),
\]
 where
\[
B_{1}A_{ij}B_{1}^{-1}=\left(\begin{array}{cc}
a_{22}^{ij} & -a_{21}^{ij}\\
-a_{12}^{ij} & a_{11}^{ij}
\end{array}\right).
\]
 Hence equations (\ref{eq:dec10.13.1}) in the component form are
reduced to the equations
\[
a_{22}^{ij}-a_{11}^{ij}+c\delta_{ij}=0,\,\,\, a_{11}^{ij}-a_{22}^{ij}+c\delta_{ij}=0,\,\,\, a_{21}^{ij}+a_{12}^{ij}=0.
\]
 Thus we find that $c=0$ and
\[
A_{ij}=\left(\begin{array}{cc}
\alpha_{ij} & \beta_{ij}\\
-\beta_{ij} & \alpha_{ij}
\end{array}\right),
\]
 where $\alpha_{ij}$ and $\beta_{ij}$ are real numbers. Due to the
commutativity of $B_{d_{1}}$ and $A_{3}$, the condition
\[
C_{0}B_{d_{1}}-B_{d_{1}}C_{0}=0
\]
provides that
\[
CB_{d_{1}}-B_{d_{1}}C=0.
\]
 In fact,
\[
CB_{d_{1}}-B_{d_{1}}C=e^{xA_{3}}C_{0}e^{-xA_{3}}B_{d_{1}}-B_{d_{1}}e^{xA_{3}}C_{0}e^{-xA_{3}}=e^{xA_{3}}(C_{0}B_{d_{1}}-B_{d_{1}}C_{0})e^{-xA_{3}}=0.
\]

Thus,
\[
C(x)=e^{xA_{3}}C_{0}e^{-xA_{3}},
\]
 where the matrices $C_0$, $A$ and $A_3$ satisfy the conditions
 \[
 C_0B_{d_1}- B_{d_1} C_0 = 0,\ \  \
A_{3}B_{d_1}-B_{d_1}A_{3}=0,\ \ \
 A_{3}C_{0}-C_{0}A_{3}\neq 0.
 \]

As an example let us consider $n=2$, and the matrices
\[
A_{3}=\left(\begin{array}{cccc}
0 & 0 & 0 & 1\\
0 & 0 & -1 & 0\\
0 & 0 & 0 & 0\\
0 & 0 & 0 & 0
\end{array}\right),\,\,\, C_{0}=\left(\begin{array}{cccc}
c_{11} & c_{12} & c_{13} & c_{14}\\
-c_{12} & c_{11} & -c_{14} & c_{13}\\
c_{31} & c_{32} & c_{33} & c_{34}\\
-c_{32} & c_{31} & -c_{34} & c_{33}
\end{array}\right).
\]
It is trivial to check that if $A_{ij}=0$ for $i\leq j$, then $A_{3}^{n}=0$,
and hence
\[
e^{xA_{3}}=\sum_{j=0}^{n}\frac{x^{j}}{j!}A_{3}^{j}.
\]
 For these matrices one obtains
\[
C_{1}\equiv C_{0}A_{3}-A_{3}C_{0}=\left(\begin{array}{cccc}
c_{32} & -c_{31} & c_{34}-c_{12} & c_{11}-c_{33}\\
c_{31} & c_{32} & c_{33}-c_{11} & c_{34}-c_{12}\\
0 & 0 & -c_{32} & c_{31}\\
0 & 0 & -c_{31} & -c_{32}
\end{array}\right)\neq0,
\]
 and
\[
C=C_{0}+xC_{1}+x^{2}\left(\begin{array}{cccc}
0 & 0 & c_{31}\,\,\,\,\,\, & c_{32}\\
0 & 0 & -c_{32}\,\,\,\,\, & c_{31}\\
0 & 0 & 0 & 0\\
0 & 0 & 0 & 0
\end{array}\right).
\]

\subsection{Classes of systems admitting three-dimensional Lie algebras}

For this case the basis of such algebras consists of the generators

\[
X_{A},\,\,\, X_{2}+X_{A_{2}},\,\,\, X_{3}+X_{A_{3}}.
\]
 The commutators of these generators are
\[
\begin{array}{c}
[X_{2}+(A_{2}\mathbf{y})\cdot\nabla,X_{3}+(A_{3}\mathbf{y})\cdot\nabla]=-2X_{3}+((A_{3}A_{2}-A_{2}A_{3})\mathbf{y})\cdot\nabla,\\
\ [X_{2}+(A_{2}\mathbf{y})\cdot\nabla,(A\mathbf{y})\cdot\nabla]=((AA_{2}-A_{2}A)\mathbf{y})\cdot\nabla,\\
\ [X_{3}+(A_{3}\mathbf{y})\cdot\nabla,(A\mathbf{y})\cdot\nabla]=((AA_{3}-A_{3}A)\mathbf{y})\cdot\nabla.
\end{array}
\]
Since the admitted Lie algebra is three-dimensional, then
\[
A_{3}(A_{2}+2E)-A_{2}A_{3}=\alpha_{1}A,\,\,\, A_{2}A-AA_{2}=\alpha_{2}A,\,\,\, A_{3}A-AA_{3}=\alpha_{3}A.
\]
 Choosing the matrix $A=B_{d_{1}}$ leads to $\alpha_{2}=0$
and $\alpha_{3}=0$. As shown earlier the condition that
the matrices $C_{0}$ and $B_{d_{1}}$ commute is sufficient for satisfying
the determining equations $CB_{d_{1}}-B_{d_{1}}C=0$ in this case.

Hence
\[
C(x)=e^{xA_{3}}C_{0}e^{-xA_{3}},
\]
 and the matrices $C_{0}$, $A_{2}$ and $A_{3}$ satisfy the conditions:

\begin{equation}
A_{3}(A_{2}+2E)-A_{2}A_{3}=\alpha_{1}B_{d_{1}},\,\,\, A_{2}B_{d_{1}}-B_{d_{1}}A_{2}=0,\,\,\, A_{3}B_{d_{1}}-B_{d_{1}}A_{3}=0,\label{eq:jan14.11.1}
\end{equation}

\begin{equation}
B_{d_{1}}C_{0}-C_{0}B_{d_{1}}=0,\,\,\, A_{2}C_{0}-C_{0}(A_{2}+4E)=0\,\,\, A_{3}C_{0}-C_{0}A_{3}\neq0,\,\,\, tr(C_{0})=0.\label{eq:jan14.11.2}
\end{equation}

\textbf{Remark}. Analyzing the first equation in (\ref{eq:jan14.11.1}).
by multiplying it by $B_{d_{1}}$, and using commutativity of the matrices
$A_{2}$ and $A_{3}$ with the matrix $B_{d_{1}}$, we obtain
\[
A_{3}B_{d_{1}}(A_{2}+2E)-A_{2}A_{3}B_{d_{1}}=-\alpha_{1}E.
\]
Then the general solution of this equation%
\footnote{See the details in Appendix.%
} is
\[
A_{3}=A_{h}-\frac{\alpha_{1}}{2}B_{d_{1}},
\]
where $A_{h}$ is the general solution of the homogeneous equation
\[
A_{h}(A_{2}+2E)-A_{2}A_{h}=0.
\]
Note that for $A_h=0$ one has that $C(x)$ is constant, then
$A_{h}\neq0$.
This implies that the matrix $A_{2}$ has an eigenvalue, say $\lambda$,
such that $\lambda+2$ is also an eigenvalue of the matrix $A_{2}$.

\subsection{Classes of systems admitting four-dimensional Lie algebras}

The basis of such Lie algebras consists of the generators

\[
X_{A},\,\,\, X_{1}+X_{A_{1}},\,\,\, X_{2}+X_{A_{2}},\,\,\, X_{3}+X_{A_{3}}.
\]
 The commutators of these generators are
\[
\begin{array}{rcl}
[X_{1}+(A_{1}\mathbf{y})\cdot\nabla,X_{2}+(A_{2}\mathbf{y})\cdot\nabla] & = & -2(X_{1}+(A_{1}\mathbf{y})\cdot\nabla)+((A_{2}A_{1}-A_{1}A_{2}+2A_{1})\mathbf{y})\cdot\nabla,\\
\ [X_{1}+(A_{1}\mathbf{y})\cdot\nabla,X_{3}+(A_{3}\mathbf{y})\cdot\nabla] & = & -(X_{2}+(A_{2}\mathbf{y})\cdot\nabla)+((A_{3}A_{1}-A_{1}A_{3}+A_{2})\mathbf{y})\cdot\nabla,\\
\ [X_{2}+(A_{2}\mathbf{y})\cdot\nabla,X_{3}+(A_{3}\mathbf{y})\cdot\nabla] & = & -2(X_{3}+(A_{3}\mathbf{y})\cdot\nabla)+((A_{3}A_{2}-A_{2}A_{3}+2A_{3})\mathbf{y})\cdot\nabla.
\end{array}
\]
\[
\begin{array}{c}
[X_{1}+(A_{1}\mathbf{y})\cdot\nabla,(A\mathbf{y})\cdot\nabla]=((AA_{1}-A_{1}A)\mathbf{y})\cdot\nabla,\\
\ [X_{2}+(A_{2}\mathbf{y})\cdot\nabla,(A\mathbf{y})\cdot\nabla]=((AA_{2}-A_{2}A)\mathbf{y})\cdot\nabla,\\
\ [X_{3}+(A_{3}\mathbf{y})\cdot\nabla,(A\mathbf{y})\cdot\nabla]=((AA_{3}-A_{3}A)\mathbf{y})\cdot\nabla.
\end{array}
\]

Since the admitted Lie algebra is four-dimensional, then
\[
AA_{1}-A_{1}A=\alpha_{1}A,\,\,\, A_{2}A-AA_{2}=\alpha_{2}A,\,\,\, A_{3}A-AA_{3}=\alpha_{3}A.
\]
 Choosing the matrix $A=B_{d_{1}}$ leads to $\alpha_{1}=0$,
$\alpha_{2}=0$ and $\alpha_{3}=0$.

Therefore
\[
C(x)=e^{xA_{3}}C_{0}e^{-xA_{3}},
\]
 and the matrices $C_{0}$, $A_{1}$, $A_{2}$ and $A_{3}$ satisfy
the conditions:
\[
B_{d_{1}}A_{1}-A_{1}B_{d_{1}}=0,\,\,\, A_{2}B_{d_{1}}-B_{d_{1}}A_{2}=0,\,\,\, A_{3}B_{d_{1}}-B_{d_{1}}A_{3}=0,
\]
\[
A_{1}A_{2}-(A_{2}+2E)A_{1}=\beta_{1}B_{d_{1}},\,\,\, A_{1}A_{3}-A_{3}A_{1}-A_{2}=\beta_{2}B_{d_{1}},\,\,\, A_{2}A_{3}-A_{3}A_{2}-2A_{3}=\beta_{3}B_{d_{1}},
\]
\[
A_{2}C_{0}-C_{0}(A_{2}+4E)=0,\,\,\, C_{0}A_{1}-A_{1}C_{0}=0,\,\,\, C_{0}B_{d_{1}}-B_{d_{1}}C_{0}=0,\,\,\, tr(C_{0})=0,\,\,\, A_{3}C_{0}-C_{0}A_{3}\neq0.
\]

As before the condition that the matrices $C_{0}$ and $B_{d_{1}}$ commute
is sufficient for satisfying the determining equations $CB_{d_{1}}-B_{d_{1}}C=0$ in this particular case.

\subsection{Systems admitting two generators $X_{A_{1}}$ and $X_{A_{2}}$ }

For completeness we also derive conditions for systems admitting a
Lie algebra with basis generators $X_{A_{1}}$ and $X_{A_{2}}$. Choosing
$A_{2}=B_{d_{1}}$, the determining equations become
\begin{equation}
CA_{1}-A_{1}C=0,\,\,\, CB_{d_{1}}-B_{d_{1}}C=0.\label{eq:dec28.13.1}
\end{equation}
 Note that
\[
A_{1}=PB_{d_{1}}P^{-1},
\]
 with some nonsingular constant matrix satisfying the condition
\[
PB_{d_{1}}-B_{d_{1}}P\neq0.
\]
 The first equation in (\ref{eq:dec28.13.1}) can be also rewritten in
the form
\[
P^{-1}(CA_{1}-A_{1}C)P=(P^{-1}CP)B_{d_{1}}-B_{d_{1}}(P^{-1}CP)=0.
\]

\section{Conclusion}
We have given a general study of the group classification of systems of linear second-order ordinary differential equations and found that irreducible linear systems (\ref{tg1}) admitting
one- two- or three-dimensional Lie algebras (\ref{eq:dec2613.1})
are equivalent to one of the following cases given in subsection (5.7), that is, $(a), (b)$ and $(c)$. This result has been stated as a Theorem in Section $5$ of the paper. The results were discussed and in two or more of the cases which had the  admitted generators of the form $X_{A}$ the study was similar to that done in the earlier studies in \cite{bk:MoyoMeleshkoOguis[2013],bk:MoyoMeleshkoOguis[2014],%
bk:MkhizeMoyoMeleshko[2014],bk:SuksernMoyoMeleshko[2014]}.
A detailed discussion on systems admitting generators of the
form $X_{A}$ has been given in Section $6$ of the paper.
We note that the further classification of
systems of second-order ordinary differential equations of
the type (\ref{tg1}), considered here, is related with the
classification of the matrix $A$, reducing it to the Jordan form.

\section{Appendix }

Here analysis of matrix equations is discussed.

\subsection{Algebraic background}

\subsubsection{Equations $CB-BC=0$}

Assume that there exists a constant real-valued $m\times m$ matrix
$B$ such that
\[
CB-BC=0.
\]
 Notice that since the matrix $B$ is real-valued, then for a complex
eigenvalue $\lambda$ of the matrix $B$ there is a conjugate eigenvalue
$\bar{\lambda}$. We call a matrix $C$ a reducible matrix if there
exists a nonsingular real-valued matrix $P$ such that
\[
PCP^{-1}=\left(\begin{array}{cc}
G_{1} & G_{2}\\
0 & G_{3}
\end{array}\right),
\]
 where $G_{1}$ and $G_{3}$ are squared matrices with $dim\, G_{3}\geq1$.

It is shown in this section that for an irreducible matrix $C$ with
$m\geq3$ the matrix $B$ can be assumed to be one of the matrices:
either $B=\alpha E$ or
\begin{equation}
B=\left(\begin{array}{cccc}
B_{0} & 0 & ... & 0\\
0 & B_{0} & ... & 0\\
... & ... & ... & ...\\
0 & 0 & ... & B_{0}
\end{array}\right),
\ \ \ B_{0}=\left(\begin{array}{cc}
\alpha & \beta\\
-\beta & \alpha
\end{array}\right),
\label{eq:dec513.0}
\end{equation}
 where $E$ is the unit matrix, $\alpha$ and $\beta$ are real numbers.

Let the matrix $B$ have at least one eigenvalue $\lambda_{1}$ such
that $\lambda_{1}$ (and its conjugate $\bar{\lambda}_{1}$ in the
case of complex $\lambda_{1}$) differs from other eigenvalues, then
by virtue of the decomposition theorem (see p.160 in \cite{bk:Shilov})
there exists a nonsingular real-valued matrix $P$ such that
\[
B=P^{-1}\left(\begin{array}{cc}
B_{1} & 0\\
0 & B_{2}
\end{array}\right)P,
\]
where the matrices $B_{1}$ and $B_{2}$ are square matrices which
have no common eigenvalues.

Using the change
\[
\tilde{\mathbf{y}}=P\mathbf{y},
\]
 and because of the equality
\[
PCP^{-1}(PBP^{-1})-(PBP^{-1})PCP^{-1}=0,
\]
 one can assume that
\[
B=\left(\begin{array}{cc}
B_{1} & 0\\
0 & B_{2}
\end{array}\right).
\]
 Let
\[
C=\left(\begin{array}{cc}
C_{1} & C_{2}\\
C_{3} & C_{4}
\end{array}\right),
\]
 then
\[
CB-BC=\left(\begin{array}{cc}
(C_{1}B_{1}-B_{1}C_{1}) & (C_{2}B_{2}-B_{1}C_{2})\\
(C_{3}B_{1}-B_{2}C_{3}) & (C_{4}B_{2}-B_{2}C_{4})
\end{array}\right)=0.
\]
 Hence, one obtains that
\[
C_{2}B_{2}-B_{1}C_{2}=0,\,\,\, C_{3}B_{1}-B_{2}C_{3}=0.
\]
 Since $B_{1}$ and $B_{2}$ have no common eigenvalues, these equations
imply that (see p.196 in \cite{bk:Gantmacher})
\[
C_{2}=0,\,\,\, C_{3}=0.
\]
 This means that a system of ordinary differential equations (\ref{tg1})
with the matrix $C $ is reducible. Thus, for irreducible systems
the matrix $B$ has only one eigenvalue, say $\lambda$
(and its conjugate $\bar{\lambda}$ in the case of complex $\lambda$).

Further analysis%
\footnote{Symbolic computer calculations were applied for this analysis.%
} gives that for irreducible to fewer dimension systems the matrix
$B$ is equivalent to a diagonal matrix. This means that the matrix
$B$ can be only of the following two types: either $B=\alpha E$ or $B=B_{d}$.
Moreover, if the matrix
$B$ is of the second type, then the matrix $C$ has a block structure
of the form ($m=2n$)
\begin{equation}
C=\left(\begin{array}{cccc}
C_{11} & C_{12} & ... & C_{1n}\\
C_{21} & C_{22} & ... & C_{2n}\\
... & ... & ... & ...\\
C_{n1} & C_{n2} & ... & C_{nn}
\end{array}\right),\label{eq:dec513.2}
\end{equation}
 where
\[
C_{ij}=\left(\begin{array}{cc}
\alpha_{ij} & \beta_{ij}\\
-\beta_{ij} & \alpha_{ij}
\end{array}\right),
\]
 $\alpha_{ij}$ and $\beta_{ij}$ are real numbers.

\subsubsection{Analysis of the matrix equation $AG-G(A+2E)=\alpha E$}

Here we consider the matrix equation
\begin{equation}
AG-G(A+2E)=\alpha E,\label{eq:dec513.5}
\end{equation}
 where the matrix $G$ and $A$ are real-valued $m\times m$ matrices
and $\alpha$ is a real number. The matrix $A$ is assumed to be given,
the matrix $G$ is unknown.

One can check that $G=-\frac{\alpha}{2}E$ is a particular solution
of this equation. Hence, the general solution of equation (\ref{eq:dec513.5})
is (see \cite{bk:Gantmacher})
\[
G=G_{h}-\frac{\alpha}{2}E,
\]
 where $G_{h}$ is the general solution of the homogeneous equation
\begin{equation}
AG-G(A+2E)=0.\label{eq:dec513.7}
\end{equation}

If the matrices $A$ and $A+2E$ have distinct sets of eigenvalues,
then the general solution of the homogeneous equation is unique (see
\cite{bk:Gantmacher}): $G_{h}=0$.

\subsection{Solution of the matrix ordinary differential equation }

The solution of the Cauchy problem of the matrix ordinary differential
equation (see p.175 in \cite{bk:Bellman})
\[
X^{\prime}=A(x)X+XB(x),\,\,\, X(0)=C_{0}
\]
 is given by the formula
\[
X=YC_{0}Z,
\]
 where
\[
Y^{\prime}=A(x)Y,\,\,\, Y(0)=E,\,\,\,\, Z^{\prime}=ZB(x),\,\,\, Z(0)=E.
\]



\end{document}